\documentclass{amsart}
\usepackage{amsmath,amsthm,amssymb,amsfonts}
\usepackage[colorlinks=true]{hyperref}
\hypersetup{urlcolor=blue, citecolor=blue}

\def\arXiv#1{   {\href{http://arxiv.org/abs/#1}
   {{\mdseries\ttfamily arXiv:#1}}}}

\def\doi#1{   {\href{http://dx.doi.org/#1}
   {{\mdseries\ttfamily DOI}}}}

\newcommand{\de}{\delta}    
  \newcommand{\ep}{\varepsilon}
\newcommand{\la}{\lambda}

\newcommand{\R}{\mathbb{R}}\newcommand{\Z}{\mathbb{Z}}

\newcommand{\M}{\mathcal{M}}

\newcommand{\pt}{\partial_t}\newcommand{\pa}{\partial}
\newcommand{\les}{{\lesssim}}
\newcommand{\beeq}{\begin{equation}}\newcommand{\eneq}{\end{equation}}

\newcommand{\supp}{\text{supp}}

\newenvironment{prf}{\noindent {\bf Proof.} }{\endprf\par}
\def \endprf{\hfill  {\vrule height6pt width6pt depth0pt}\medskip}

\numberwithin{equation}{section}

\def\Hr{ H_{{\rm rad}}}

\def\<{\langle}             \def\>{\rangle}
\def\({\left(}                 \def\){\right)}

\newtheorem{thm}{Theorem}[section]

\newtheorem{coro}[thm]{Corollary}
\newtheorem{lem}[thm]{Lemma}

\theoremstyle{remark}
\newtheorem{rem}{Remark}[section]

\theoremstyle{definition}
\title[Weighted fractional chain rule and
nonlinear wave equations]
{
Weighted fractional chain rule and
nonlinear wave equations with minimal regularity
}

\author{Kunio Hidano}
\address{Department of Mathematics\\
Faculty of Education\\
Mie University\\
1577 Kurima-machiya-cho, Tsu, Mie 514-8507, JAPAN}
\email{hidano@edu.mie-u.ac.jp}

\author{Jin-Cheng Jiang}
\address{Department of Mathematics, 
         National Tsing Hua University, 
          Hsinchu, Taiwan 30013, R.O.C.}
\email{jcjiang@math.nthu.edu.tw}

\author{Sanghyuk Lee}
\address{Department of Mathematical Sciences\\ Seoul National University\\ Seoul 151-747, Republic of Korea}
\email{shklee@snu.ac.kr}

\author{Chengbo Wang}\address{School of Mathematical Sciences\\                Zhejiang University\\                Hangzhou 310027, China}\email{wangcbo@gmail.com}
\urladdr{http://www.math.zju.edu.cn/wang}

\date{\today}
\dedicatory{} \commby{}

\begin{document}

\begin{abstract} We consider the local well-posedness for 3-D quadratic semi-linear wave equations with radial data: 
\[  \begin{aligned}
&\qquad \qquad \Box u = a |\pt u|^2+b|\nabla_x u|^2, \\
& u(0,x)=u_0(x)\in H^{s}_{\mathrm{rad}}, \  \pt u(0,x)=u_1(x)\in H^{s-1}_{\mathrm{rad}}.
\end{aligned}\] 
It has been known that the problem is well-posed for $s\ge 2$ and ill-posed for $s<3/2$.  In this paper,  we prove unconditional well-posedness up to the scaling invariant regularity, that is to say,  for $s>3/2$ and thus fill the gap which was left open for many years.  For the purpose, we also obtain a weighted fractional chain rule, which is of independent interest. Our method here also works for a class of nonlinear wave equations with general power type nonlinearities which contain  the space-time derivatives of the unknown functions. In particular, we  prove the Glassey conjecture in the radial case, with minimal regularity assumption.
\end{abstract}

\keywords{Glassey conjecture, fractional chain rule, nonlinear wave equations, generalized Strichartz estimates, unconditional uniqueness}

\subjclass[2010]{35L70, 35L15, 42B25, 42B37}

\maketitle 
%\tableofcontents

\section{Introduction}
 In this paper, we are interested in the problem of local well-posedness with minimal regularity for 3-D quadratic semi-linear wave equations
with radial data, 
\begin{equation}
\begin{aligned}\label{eq-nlw}
&\qquad \qquad \Box u = a |\pt u|^2+b|\nabla_x u|^2, \\
& u(0,x)=u_0(x)\in H^{s}_{\mathrm{rad}}, \  \pt u(0,x)=u_1(x)\in H^{s-1}_{\mathrm{rad}},
\end{aligned}
\end{equation} 
where $(t,x)\in \R_+\times\R^{3}$, $(a,b)\in\R^2\backslash \{(0,0)\}$, $\Box =\pt^2-\Delta$, $\pa=(\pt, \nabla_x)$.
Here $H^s_{\rm{rad}}$  denotes 
the space of spherically symmetric functions lying in the usual Sobolev space $H^s$. 
Notice that
$s_c=3/2$ is the regularity for the problem to be scale invariant in the homogeneous Sobolev space $\dot H^{s}$, and it is well known that the problem is ill-posed in $H^s$ and $\Hr^s$ for $s<s_c$, see Lindblad \cite{Ld93}.
The equation \eqref{eq-nlw} serves as a simplified model for many equations which are important  in mathematical physics, such as wave maps and Einstein equations. The problem of local well-posedness with minimal regularity assumption has been extensively investigated and well-understood.

\subsection{Quadratic semi-linear wave equation} In \cite{KlMa93}, by using 
bilinear estimates, together with standard energy-type estimates,
Klainerman and Machedon proved the local well-posedness for the problem \eqref{eq-nlw} with $a+b=0$ (thus, satisfying the null condition) in $H^2$,
which was later improved to $H^s$ for $s>s_c$, see \cite{KM95},  \cite{KM96}, \cite{KS97} and references therein. 
 Moreover, 
 even without the null condition,
their bilinear estimate holds
if initial data and forcing terms are
radially symmetric in the spatial variables.
It yields local well-posedness in $\Hr^2$ for the problem \eqref{eq-nlw},
  and they conjectured that the problem is ill posed in $H^2$ in general.
  Later, 
Ponce-Sideris \cite{PS} proved local well-posedness in $H^s$ for any $s>2$, as well as the radial case $\Hr^2$, by applying the Strichartz estimates.
See also Tataru \cite{Tataru99}, Smith-Tataru \cite{SmTa05} and references  therein for related results on local well-posedness  for  the same problem as well as the quasilinear wave equation $\Box_{g(u)} u=Q(\pa u, \pa u)$, in general spatial dimensions. 
In \cite{Ld93}, \cite{Ld96}, Lindblad proved that the problem is {generally} ill-posed in $H^s$ for any $s\le 2$, which verified the conjecture of Klainerman-Machedon \cite{KlMa93}.
Thus, in general cases  the local well-posedness problem looks almost completely understood in terms of regularity of the initial data. However, if we focus on the radial data, the  current state is not so satisfactory.  To be precise, the results of Lindblad \cite{Ld93}, \cite{Ld96} do not apply to the radial solutions, so we only know that the problem is ill-posed in $\Hr^s$ with $s<s_c=3/2$ (see {\it Remark} \ref{rem1}). On the side of well-posedness,  from \cite{KlMa93}  the problem \eqref{eq-nlw} is locally well-posed in $\Hr^2$.  This means, in contrast to the general cases, that  for the radial case  
there is still a gap between ill-posedness and  well-posedness results,  precisely 
for \eqref{eq-nlw} in $\Hr^s$ with $s\in [3/2, 2)$.

In the following which is our main result  we prove that $s>s_{c}$ is sufficient for \eqref{eq-nlw} to be locally well-posed  so that we fill the gap except the case $s=s_c$. Moreover, we show that the local solutions could be extended to almost global solutions in $\Hr^s$, when the initial data are small.  Since the result for $s\ge 2$ was already proved in the aforementioned previous works, we restrict ourselves to the case $s\in (s_c, 2)=(3/2, 2)$.

\begin{thm}
\label{thm-main0}
The problem  \eqref{eq-nlw} is unconditionally well-posed in $\Hr^s$ with $s\in (3/2,2)$. More precisely, for any radial data $(u_0,u_1)\in \Hr^s\times \Hr^{s-1}$,  there exists $T\in (0,\infty)$, which depends only on the norm of $(u_0, u_1)$, such that the problem \eqref{eq-nlw} admits a unique, radially symmetric solution 
$u\in L^\infty([0, T]; H^s(\R^3))\cap Lip( [0, T]; H^{s-1}(\R^3))$.
Moreover,  $u\in C([0, T]; H^s)\cap C^1( [0, T]; H^{s-1})$,
$$r^{-1/2+\de}\<r\>^{-\de} \pa u\in L^2([0, T]\times \R^{3}),\  2\de=s-\frac{3}{2}\ ,$$
 and there exists $c>0$, such that we can choose $T=T_\ep$, 
where 
\beeq\label{eq-life}T_\ep=
\left\{
\begin{array}{ ll}
  \exp(c\ep^{-1}) \ ,&   \ep<1\ ,\\
c\ep^{-\frac{1}{s-3/2}}\ ,  &   \ep\ge 1\ ,
\end{array}
\right.
\eneq
 with 
$$
\ep=
  \|\pa u(0)\|_{\dot H^{s-1}}+ \|\pa u(0)\|_{\dot H^{s-1}}^{1/2}
\|\pa u(0)\|_{\dot H^{2-s}}^{1/2}.
$$
\end{thm}

We now make several remarks concerning  Theorem \ref{thm-main0}. 
\begin{rem} 
It is natural to expect  that solutions exist for longer time as the initial data get smaller but it is another matter to obtain precise quantification of  the lifespan depending on the size of the initial data.  For the problem \eqref{eq-nlw} with generic compactly supported smooth data of size $\ep$, it is known that the lifespan $T_\ep$ of the classical solution satisfies
$$\ln T_\ep\simeq \ep^{-1},$$
except the case when the null conditions are satisfied (in which case we have small data global existence), see
John \cite{Jo81}, John-Klainerman \cite{JK84}, 
Rammaha \cite{Ram97},
Zhou \cite{Zh01}.
For $\Hr^2$ data 
Hidano-Yokoyama \cite{HiYo06} proved almost global existence, i.e., $T_\ep\ge \exp (c\ep^{-1})$, of the Klainerman-Machedon radial solutions for small $\Hr^2$ data by using a variant of the KSS estimate \cite{KSS}. See Hidano-Wang-Yokoyama  \cite{HWY1} for {a related result}  in the quasilinear case. See also Sterbenz \cite{Ster07} and references therein for small data global well-posedness with low regularity, for \eqref{eq-nlw} with spatial dimension $n\ge 4$.
\end{rem}

\begin{rem}\label{rem1}
The regularity assumption in Theorem \ref{thm-main0}, $s>s_c=3/2$, is optimal, in the sense that the same conclusion fails for $s=s_c$ in general. More precisely, for \eqref{eq-nlw} with $a=1$, there exists a family of compactly supported, smooth, radially symmetric functions $g_\la$, which are uniformly bounded in $\Hr^{1/2}$, such that the lifespan of the corresponding solution $u_\la$ in $\Hr^s$, with data $u(0)=0$, $u_t(0)=g_\la$, goes to zero as $\la$ goes to $\infty$. See the proof of Fang-Wang \cite[Theorem 1.2]{FW1} together with an explicit ODE solution $u_\la(x,t)=-\ln (1-t/\la)$ instead of (3.3) there.
\end{rem}

\begin{rem}
The unconditional uniqueness for 
the problem  \eqref{eq-nlw} has been obtained in $H^s$ with $s>2$ in Planchon
\cite{Plan03}. See also Zhou
\cite{Zh00} for a related $H^2$ result under the null condition.
\end{rem}

\subsection{Weighted fractional chain rule}   
The well-posedness results in  \cite{KlMa93}, \cite{PS} are based on bilinear and Strichartz estimates. As is clear, such approach doesn't seem to work for our problems. Instead, for the proof of Theorem \ref{thm-main0},  we make  use  of the Morawetz type estimates for the wave equations (also known as local energy estimates and KSS type estimates in the literature) and  a weighted fractional chain rule which is   
of independent interest. As far as we are aware, this type of chain rule hasn't appeared before.
 As Morawetz type estimates could be proved by multiplier methods (for example, see \cite{Ster05, MetSo06}) our proof of Theorem \ref{thm-main0} provides a {physical space approach}  for the nonlinear wave equations with low regularity. We believe that our approach will be useful for studies in related nonlinear problems.

As for the weight functions, we use Muckenhoupt's  $A_p$ class.   Let us recall the definition of  $A_p$ weights: 
$$w\in A_1\Leftrightarrow \mathcal{M}w(x)\le Cw(x),  \  a.e.\ x \in\R^n\ ,$$
and,  for $1<p<\infty$,  $w\in A_p$ if and only if $w\ge 0$  and 
$$ \left(\int_Q w(x)dx\right)\left(\int_Q w^{1-p'}(x)dx\right)^{p-1}\le C|Q|^p, \  \   \forall \textrm{ cubes }Q\,.$$
Here  $\M$ is the Hardy-Littlewood Maximal function which is defined  by $\M w(x)=\sup_{r>0}r^{-n}\int_{B_r(x)} |w(y)|dy$.

\begin{thm}[Weighted fractional chain rule]\label{thm-wLeib}
Let $s\in (0,1)$, $q, q_{1}, q_{2}\in (1,\infty)$ with $\frac 1q=\frac{1}{q_{1}}+\frac{1}{q_{2}}$.
Assume $F:\R^k\rightarrow \R^l$ is a $C^1$ map, satisfying $F(0)=0$ and 
\beeq\label{eq-wLeib-assu}|F'(\tau v+(1-\tau)w)|\le \mu(\tau)|G(v)+G(w)|,\eneq
with $G\ge 0$ and $\mu\in L^1([0,1])$.
If
$(w_1 w_{2})^{q}\in A_{q}$, $w_{1}^{q_{1}}\in A_{q_{1}}$,
$w_{2}^{q_2}\in A_{q_{2}}$, then 
\beeq\label{eq-wLeib1}\|w_1w_{2} D^s F(u)\|_{L^q}\les \|w_1 D^s u\|_{L^{q_{1}}} \|w_{2} G(u)\|_{L^{q_{2}}}\ ,\eneq
where $D=\sqrt{-\Delta}$.
In addition, when $q_{2}=\infty$ and $q\in (1,\infty)$, if
$w_1^{q}, (w_{1}w_2)^{q}\in A_q$ and $w_2^{-1}\in A_1$, we have
\beeq\label{eq-wLeib}\|w_1 w_{2} D^s F(u)\|_{L^q}\les \|w_1 D^s u\|_{L^q} \|w_2 G(u)\|_{L^\infty}\ .\eneq
\end{thm}

The fractional chain rule with $w_1=w_2\equiv 1$(without weights) was proven by
Christ-Weinstein \cite{CW91}, Staffilani \cite{Sta97}, Kato \cite{Ka95}.
See also Taylor  \cite[Chapter 2 Proposition 5.1]{Tay}.
Our proof of Theorem \ref{thm-wLeib} basically follows the same line of agrument as in Taylor \cite{Tay} but we  also need new ingredients such as the weighted estimates for Calder\'on-Zygmund, Littlewood-Paley operators and weighted vector valued inequality for  Hardy-Littlewood maximal functions (see 
Andersen-John \cite{AJ80}). 
We note the further potential of
the fractional chain rule in Theorem \ref{thm-wLeib} and
the characterization of our weight
$w(x)=r^{-1+2\delta}\<r\>^{-2\delta-2\delta_1}$
in terms of the Muckenhoupt $A_1$ class
(see Lemma \ref{thm-Ap} below).
They also play an important role in the study
of well-posedness of the Cauchy problem for
the nonlinear half wave equation and the nonlinear elastic wave equation with low-regularity data. See Hidano-Wang
\cite{HW17}, Hidano-Zha \cite{HZ17}.

\subsection{General power type nonlinearities} We also apply our method to a class of nonlinear wave equations with general power type nonlinearities in which  the space-time derivatives of the unknown functions also appear. More precisely, let $n\ge 2$, $p>1$, and consider the following nonlinear wave equation, with $(t,x)\in \R_+\times\R^{n}$,
\begin{equation}
\begin{aligned}
\label{a-pde}
&\qquad\qquad \quad\Box u = a |\pt u|^p+b|\nabla u|^p, \\
&u(0,x)=u_0(x)\in H^{s}_{\mathrm{rad}}, \  \pt u(0,x)=u_1(x)\in H^{s-1}_{\mathrm{rad}}.
\end{aligned}
\end{equation} 
For this problem, $s_c=(n+2)/2-1/(p-1)$ is the scaling-critical Sobolev regularity, which is a lower bound for the problem to be locally well-posed in $H^s$, and also $\Hr^s$. In addition, $s_l=(n+5)/4$ is a lower bound for the problem to be locally well-posed in $H^s$.
See Fang-Wang \cite{FW1} for these ill-posed results when $n\ge 3$\footnote{
In Fang-Wang \cite{FW1},
the theorems were stated for equations of the form $\Box u=|u_t|^{p-1} u_t$, however, as 
 is clear from the proof, it also applies to \eqref{eq-nlw} with $a>0$ and $b\ge 0$.}.
If the nonlinearity is sufficiently smooth, the problem is locally well-posed in $H^s$, for $s\ge s_c$ and $s>\max(s_l, (n+1)/2)$. This can be shown by using Strichartz estimates, see, e.g., Fang-Wang \cite{FW1}. Moreover, when the initial data are radial,  {improved Strichartz type estimates under radial symmetry} can be further exploited to prove local well-posedness in $\Hr^s$ for  $s\ge \max(s_c, (n+1)/2)$ and $s>3/2$ (see, Fang-Wang \cite{FW2}, \cite{FW13}).

Concerning long time existence with small data, the problem is related to the so-called Glassey conjecture (see Glassey's MR review of
the paper of Sideris \cite{Si83}),
which states that the critical power for \eqref{a-pde} to admit global small solutions is 
$$p_{c}(n)\equiv 1+\frac{2}{n-1}\ .$$
The global existence for $p>p_c(n)$ has been verified {for dimension two and three (Hidano-Tsutaya 
 \cite{HidTsu95} and
Tzvetkov \cite{Tz98}),
 as well as {radial cases in  higher dimensions}
(Hidano-Wang-Yokoyama
\cite{HWY2})}.
Notice however that the method which relies on Strichartz estimates yields small data global existence in $H^{s_c}$ for $p>\max((n+3)/(n-1), 3)$.
For the sub-critical case $1<p\le p_{c}(n)$, when $(a,b)=(1,0)$,
it is known
from Zhou \cite{Zh01}
 that there are no global solutions in general and an upper bound of the lifespan is given by
$$T_{\ep}\le \left\{
\begin{array}{ ll}
  \exp(C\ep^{-(p-1)}) \ ,&    p=p_c\ ,\\
C \ep^{\frac{2(p-1)}{(n-1)(p-1)-2}}\ ,  &   1<p<p_c\ ,
\end{array}
\right.
$$
for generic,  compactly supported, smooth data of size $\ep$.

Concerning long time well-posedness with small $\Hr^s$ data, it is known from \cite{HWY2} that $s=2$ is sufficient for existence. 
Moreover, the smallness of the initial data was measured in certain ``multiplicative form'', which strongly suggests that, the minimal regularity for the problem {should be} given by
\beeq\label{so}s_{o}\equiv \max\left(\frac{3}{2}, s_c\right)
=\max\left(\frac{3}{2}, \frac{n+2}2-\frac{1}{p-1}\right)
\ .\eneq
 Notice that for $p>1$, $s_c> 3/2$ is equivalent to $p>p_c$.

By making use of the similar approach which is taken for \eqref{eq-nlw}, we prove that the problem
\eqref{a-pde} is 
 locally well-posed in $\Hr^s$ for any $s\in (s_{o}, 2)$, when $1<p<1+2/(n-2)$. Furthermore, {when $p_c<p<1+2/(n-2)$,} for any $s\in (s_c, 2)$, it admits global solution for small data in $\Hr^s$, which gives an affirmative answer for the natural regularity problem raised in \cite{HWY2}. 

The following is our result concerning the problem \eqref{a-pde}. 

\begin{thm}
\label{thm-main1}
Let $n\ge 2$ and $p\in (1,1+2/(n-2))$ {\rm (}when $n=2$, it {means}  $p\in (1,\infty)$~{\rm )}.
The problem  \eqref{a-pde} 
 is unconditionally well-posed in $\Hr^s$ with $s\in (s_o,2)$, where $s_o$ is given by \eqref{so}. More precisely, for any radial data $(u_0,u_1)\in \Hr^s\times \Hr^{s-1}$,  there exists $T\in (0,\infty]$, which depends only on the norm of $(u_0, u_1)$, such that the problem \eqref{a-pde} admits a unique, radially symmetric solution 
$u\in L^\infty([0, T); \Hr^s)\cap Lip( [0, T); \Hr^{s-1})$. 
Moreover,  $u\in C([0, T); H^s)\cap C^1( [0, T); H^{s-1})$,
$$r^{-\frac 12+\de} 
\<r\>^{-\de-\de_{1}}
\pa u\in L^2([0, T)\times \R^{n}),
$$
where  $\de\in (0, 1/2)$ and $\de_1\in [-\de, (n-1)/2)$
are defined by
\beeq\label{eq-de}\de=\left\{
\begin{array}{ ll}
\frac{p-1}{2}(s-s_c) \ ,     &   p\ge p_c\ , \\
\frac{p-1}{2}(\frac{3}{2}-s_c) \ ,     &   p< p_c\ ,
\end{array}
\right.
\quad 
\de_1=\left\{
\begin{array}{ ll}
\frac{n-1}4 \min(1,p-p_c)\ ,& p>p_c\ ,\\
0 \ ,     &   p= p_c\ ,\\
-\de \ ,     &   p< p_c\ . 
\end{array}
\right.
\eneq
 and  there exist $c, \ep_0>0$, such that we can choose $T=T_\ep$, 
where 
\beeq\label{life}T_\ep=
\left\{
\begin{array}{ ll}
c\ep^{-\frac{1}{s-s_c}}\ ,  &   \ep\ge \ep_0, p\ge p_c\ ,\\
  \exp(c\ep^{-(p-1)}) \ ,&   \ep<\ep_0, p=p_c\ ,\\
\infty\ ,&    \ep<\ep_0, p>p_c\ ,\\
c\ep^{-\frac{1}{3/2-s_c}}=c\ep^{\frac{2(p-1)}{(n-1)(p-1)-2}}
\ ,  &   p<p_c\ ,
\end{array}
\right.
\eneq
 with 
$$
\ep=
\left\{
\begin{array}{ ll}
  \|\pa u(0)\|_{H^{s-1}}\ ,& p>p_c\ ,\\
  \|\pa u(0)\|_{\dot H^{s-1}}+ \|\pa u(0)\|_{\dot H^{s-1}}^{1/2}
\|\pa u(0)\|_{\dot H^{2-s}}^{1/2} \ ,&p=p_c\ ,\\
\|\pa u(0)\|_{\dot H^{s-1}}^{1/2}
\|\pa u(0)\|_{\dot H^{2-s}}^{1/2}\ ,& p<p_c\ .
\end{array}
\right.
$$
\end{thm}
\begin{rem}
For the two dimensional critical case $p=p_c(2)=3$, Fang-Wang \cite{FW13} proved the same result for general data in $H^{s}$ with $s>s_{c}$, by assuming additional angular regularity\footnote{As for the meaning of additional angular regularity,
see Lemma \ref{thm-trace} below. There, we use the Sobolev space
on the unit sphere} of order
 bigger than $1/2$.
\end{rem}

As we have pointed out in {\it Remark} \ref{rem1}, Theorem \ref{thm-main1} for $p\ge p_c$ fails to be true for $s=s_c$, in the sense that there exists
a family of data such that the lifespans tend to zero while the data remain bounded in $\Hr^{s_c}$. Nevertheless, it does not exclude the possibility of small data global existence in the critical space $\Hr^{s_c}$, for $p>p_c$. Thus, one may naturally  ask whether it is possible to obtain small data global existence for $s=s_c$ under the assumption $p>p_c$.  It turns out that this is true when the spatial dimension is two, while the higher dimensional cases remain open.

By exploiting generalized Strichartz estimates of 
 Fang-Wang \cite{FW2},
Smith-Sogge-Wang \cite{SSW12},   we obtain the following:
\begin{thm}
\label{thm-1}
Let $n=2$ and $p>5$. Then there exists a small constant $\ep_0>0$, such that the Cauchy problem \eqref{a-pde} has a unique global solution satisfying $u\in C([0,\infty); H^{s_c})\cap C^1([0,\infty); H^{s_c-1})$ and $\pa u\in L^{p-1}_t L^\infty_x$, whenever the initial data
$(u_0,u_1)\in H^{s_c}\times H^{s_c-1}$ with
\begin{equation}
\label{eq-small}
\|(u_0,u_1)\|_{
\dot H^{s_c}\times \dot H^{s_c-1}
}=\ep \le  \ep_0\ .
\end{equation} In addition, there exists a consant $C>0$, such that the solution satisfies
$$\|\pa u\|_{L^\infty_t \dot H^{s_c-1}\cap L^{p-1}_t L^\infty_x}\le C \ep, \   \|\pa u\|_{L^\infty_t L^2_x}\le C 
\|(u_0,u_1)\|_{
\dot H^{1}\times L^2}\, . 
$$ Moreover,  when $p>3$ and the initial data are radial, the same result  remains valid 
for the radial solutions.
\end{thm}
As we have recalled, $s_c$ and $s_l$ are believed to be lower bounds for the problem to be locally well-posed in $H^s$, which has been partially verified in 
 Lindblad \cite{Ld93}, \cite{Ld96},
Fang-Wang \cite{FW1}, for $n\ge 3$. The necessity of $s\ge \max(s_c, s_l)$ for $n=2$ is more delicate, 
but nevertheless is still true, see
Liu-Wang \cite{LW17}.
Since $s_l > s_c$ for $p < 5$,
we see that the condition $p\ge 5$ is necessary in general
for  \eqref{a-pde} to be well-posed in $H^{s_c} \times H^{s_c-1}$.

\subsubsection*{Concerning the paper} Our paper is organized as follows. In the next section, we
collect various basic estimates {which we need}, including trace estimates, Morawetz type estimates and Strichartz type estimates. In particular, we prove a weighted fractional chain rule, Theorem \ref{thm-wLeib}.
Then in Section 
\ref{Sec:3}, we {prove} local well-posedness for
\eqref{eq-nlw} and \eqref{a-pde}, 
Theorem \ref{thm-main0} and Theorem \ref{thm-main1}.
In the last section, we  prove Theorem \ref{thm-1},  the
small data global existence with critical regularity for dimension two.

\section{Preliminaries}
In this section, we collect various basic estimates to be used for proofs of theorems. All of these estimates are well known, except a weighted fractional chain rule,
Theorem \ref{thm-wLeib}.

\subsection{Trace estimates: spatial decay}
At first, let us record the trace estimates, which will provide spatial decay for functions, 
see, e.g., (1.3), (1.7) in Fang-Wang \cite{FW11}.
\begin{lem}[Trace estimates]\label{thm-trace}
Let $n\ge 2$ and $1/2< s<n/2$. Then we have
\beeq\label{eq-trace0}
\|r^{n/2-s} f\|_{L^\infty_r H^{s-1/2}_\omega}\les \|f\|_{\dot H^s}, \ 
\|r^{(n-1)/2} f\|_{L^\infty_r L^{2}_\omega}\les \|f\|_{\dot B^{1/2}_{2,1}},
\eneq
 for any $f\in C_0^\infty (\R^n)$.
Here $\dot B^{s}_{p,q}$ is the homogeneous Besov space {\rm (}see, e.g., \cite{FW2} for the definition{\rm )} and $H^s_\omega=(1-\Delta_\omega)^{-s/2}L^2_\omega$ is the Sobolev space on the unit sphere.
In particular, when $u$ is spatially radial, then $|\pa u|\les \|\pa u(t,r\cdot)\|_{L^2_\omega}$ and 
\beeq\label{eq-trace1}
\|r^{n/2-s} \pa u\|_{L^\infty_{t,x}}\les \|\pa u\|_{L^\infty_t \dot H^s}, \ 
\|r^{(n-1)/2} \pa u\|_{L^\infty_{t,x}}^2\les \|\pa u\|_{L^\infty_t \dot H^s}\|\pa u\|_{L^\infty_t \dot H^{1-s}}.
\eneq
\end{lem}

\subsection{Space-time estimates}
We will need to exploit the following space-time estimates for solutions to the linear wave equations: Morawetz type (local energy) estimates, as well as the generalized Strichartz estimates.

At first, we record the required Morawetz type estimates for the operator $\Box$. 
\begin{lem}\label{thm-LE0}
Let $n\ge 2$. Then
for any $\de\in (0, 1/2]$ and $\de_1>0$, there exists $C>0$, such that we have
\beeq\label{eq-kss3}
T^{-\de}
\|r^{-\frac 12+\de} \pa u\|_{L^2_T L^2_x}
+\|\pa u\|_{L^\infty_T L^2_x}
\le C\|\pa u(0)\|_{L^2_x}+C T^{\de}
\|r^{\frac 12-\de} \Box u\|_{L^2_T L^2_x}\ ,\eneq
\beeq
\begin{split}
 A_T^{-1}
\|r^{-\frac 1 2 +\de}\<r\>^{-\de} \pa u\|_{L^2_T L^{2}_{x}}
&+\|\pa u\|_{L^\infty_T L^2_x}
\\
& \le   C \|\pa u(0)\|_{L^2_x}+CA_T\|r^{\frac 1 2-\de}\<r\>^{\de}\Box u\|_{L^2_T L^2_x}\ ,
 \end{split}
 \label{eq-kss2}
\eneq
\beeq\label{eq-kss1}\begin{split}
B_T^{-1}\|r^{-\frac 12+\de}\<r\>^{-\de-\de_1} \pa u\|_{L^2_ T L^2_x}&+\|\pa u\|_{L^\infty_T L^2_x}\\&
\le C\|\pa u(0)\|_{L^2_x}+C B_T \|r^{\frac 12-\de}\<r\>^{\de+\de_1}\Box u\|_{L^2_TL^2_x}\ , \end{split}\eneq
for any $T\in (0,\infty)$,
where
$A_T=\min(T^\de, \sqrt{\ln(2+T)})$,
$B_T=\min(T^\de, 1)$,
 $L^{q}_{T}=L^{q}([0,T))$ for the variable $t$
and $C$ is independent of $T>0$.
\end{lem}
The  above estimates were formulated and proved in
 \cite[Lemma 3.2]{HWY2}
 for $n\ge 3$,
by multiplier method,
see also \cite{
MetSo06, HWY1, W15}.
 We remark that the estimates with
 $n\ge 2$ is also implied by the 
local energy estimates in Metcalfe-Tataru \cite[Theorem 1]{MeTa12MA}, and a well-known argument of Keel-Smith-Sogge \cite{KSS}.
Morawetz type estimates and local energy estimates have rich history and a large body of related literature. We refer the interested readers to \cite{MeTa12MA, LMSTW} for more exhaustive treatment of such estimates.

For the problem with dimension two and critical regularity, we will also use the following generalized Strichartz estimates due to Smith-Sogge-Wang \cite{SSW12}, along with the previous radial estimates in Fang-Wang \cite{FW2}.
See also \cite{JWY12} for similar results in higher dimensions.
 For Strichartz estimates, see \cite{FW2} and references therein. 
\begin{lem}[Generalized Strichartz estimates]\label{thm-SSW}
Let $n=2$, $q\in (2, 4]$, and $s=1-1/q$.  
Then we have the following inequality
\beeq\label{eq-SSW}
    \|\pa u\|_{L^q ([0,\infty);     L^\infty_r L^{2}_\omega(\R^2))}\les \|\pa u(0)\|_{\dot H^{s}_x}+\|\Box u\|_{L^1_t \dot H^{s}_x}\ .
\eneq  
In addition, we have the classical Strichartz estimates for $q\in (4,\infty)$,
\beeq\label{eq-Stri}
    \|\pa u\|_{L^q ([0,\infty); L^\infty_x (\R^2))}\les \|\pa u(0)\|_{\dot H^{s}_x}+\|\Box u\|_{L^1_t \dot H^{s}_x}\ .
\eneq 
In particular, when $u$ is {\rm (}spatially{\rm )} radial, we have
\beeq\label{eq-SSW2}
    \|\pa u\|_{
    L^q ([0,\infty); L^\infty_x (\R^2))
    }\les \|\pa u(0)\|_{\dot H^{s}_x}+\|\Box u\|_{L^1_t \dot H^{s}_x},\ q>2\ .
\eneq  
\end{lem}

As mentioned in Lemma \ref{thm-trace}, we have
$|\partial u(t,x)|
\leq
C\|\partial u(t,r\cdot)\|_{L^2_\omega}$,
when $u$ is spatially radial.
Combining \eqref{eq-SSW} with it,
we easily obtain \eqref{eq-SSW2},
which means that
the $2$-D radial estimate due to Fang-Wang (see Theorem 4 in 
\cite{FW2})
remains true for the first derivatives of radial solutions.

\subsection{Weighted fractional chain rule}
To apply the Morawetz type estimates for the nonlinear problem,
we are naturally required to introduce the weighted fractional chain rule, Theorem \ref{thm-wLeib}.
Before proving it in Section \ref{sec-3}, we firstly discuss the implication which will be useful for our problem.
 
As a direct corollary of Theorem \ref{thm-wLeib}, with $w_2=w_1^{-2}=w^{-1}$ and the fact that $$w\in A_1\Rightarrow w\in A_2
\Rightarrow w^{-1}\in A_2
\ ,$$ we obtain the following:
\begin{coro}\label{thm-wLeib2}
Let $w\in A_1$ and $s\in (0,1)$. Under the same assumption on $F(u)$ as in Theorem \ref{thm-wLeib}, we have
\beeq\label{eq-wLeib2}\|w^{-1/2} D^s F(u)\|_{L^2_x}\les \|w^{1/2} D^s u\|_{L^2_x} \|w^{-1} G(u)\|_{L^\infty_x}\ .\eneq
\end{coro}

The actual weight function we will choose is of the form
 $w(x)=r^{-1+2\de}\<r\>^{-2\de-2\de_1}$.
 \begin{lem}\label{thm-Ap}
Let 
 $w(x)=r^{-1+2\de}\<r\>^{-2\de-2\de_1}$, with
 $0\le 1-2\de\le 1+2\de_1<n$.
Then $w\in A_1(\R^n)$.
\end{lem}
\begin{prf} 
Though proof of  this lemma is rather elementary, for the sake of completeness we provide a proof.
It amounts to proving that for any $r>0$, and almost every $x\in \R^n$,
\beeq  \label{eq:a1}\int_{B_r(x)} |y|^{-1+2\de}\<y\>^{-2\de-2\de_1} dy\le C r^n |x|^{-1+2\de}\<x\>^{-2\de-2\de_1} \ .
\eneq

We deal with two cases separately. First, if $|x|\le 1$, then as $\de+\de_1\ge 0$, we have
$$|y|^{-1+2\de}\<y\>^{-2\de-2\de_1}\le |y|^{-1+2\de}\in A_1$$
provided that {$1-2\de\in [0,n)$}
(recall that $|x|^a\in A_1$ iff $a\in (-n,0]$, 
Grafakos \cite[Example 7.1.7, page 506]{Gra14}). Thus,
\[ r^{-n}\int_{B_r(x)} |y|^{-1+2\de}\<y\>^{-2\de-2\de_1} dy\le C  |x|^{-1+2\de}
\le C |x|^{-1+2\de}
\<x\>^{-2\de-2\de_1} \ , \forall r>0 .\]

For the case $|x|\ge 1$, 
recall that $|x|^{-1-2\de_1}\in A_1$ if {$1+2\de_1\in [0,n)$}. So, we have
\beeq\label{eq-A11} r^{-n}\int_{B_r(x)} |y|^{-1-2\de_1} dy\le C  |x|^{-1-2\de_1}
\ , \forall r>0 .\eneq
If
$r<|x|/2$, we have $|y|\ge |x|-r\ge |x|/2$ and so $|y|\simeq \<y\>\simeq |x|$ for $y\in B_r(x)$. Hence, 
$$ r^{-n}\int_{B_r(x)} |y|^{-1+2\de}\<y\>^{-2\de-2\de_1} dy\le 
C r^{-n}\int_{B_r(x)} |x|^{-1-2\de_1} dy\le 
 C |x|^{-1+2\de}
\<x\>^{-2\de-2\de_1} \ .$$
On the other hand, if $r\ge |x|/2$, by \eqref{eq-A11},  it follows that 
\begin{eqnarray*} 
\int_{B_r(x)} |y|^{-1+2\de}\<y\>^{-2\de-2\de_1} dy
&\le &
\int_{B_1} |y|^{-1+2\de} dy+
\int_{B_r(x)\backslash B_1} |y|^{-1-2\de_1} dy\\
&\le& 
C+C r^{n} |x|^{-1-2\de_1} \\
&\le& 
C r^{n} |x|^{-1-2\de_1} \ .
\end{eqnarray*}
This gives the desired \eqref{eq:a1} and completes the proof. 
\end{prf}

\subsection{Proof of Theorem \ref{thm-wLeib}}\label{sec-3}

To begin with, we recall that if $T$ is a strong Calder\'on-Zygmund operator, then
\beeq\|T(f)\|_{L^p(wdx)}\le C\|f\|_{L^p(wdx)}, \  w\in A_p,  \  p\in (1,\infty)\ .\eneq   (See, e.g., Muscalu-Schlag \cite[Theorem 7.21, page 191]{MuSch13}.) 
Using this and  the argument in \cite[Section 8.2]{MuSch13}   it is easy to see that the weighted Littlewood-Paley square-function estimate 
\beeq\label{eq-LPw}\|w S_j f\|_{L^p \ell^2_j}\simeq \|w f\|_{L^p}, \  w^{p}\in A_p, \  f\in L^p(w^p dx), \ p\in (1,\infty)\eneq
holds where $S_j=\phi_j*$ is the standard Littlewood-Paley operator while
$\phi_j(x)=2^{jn}\phi(2^j x)$,
 $\supp\ \hat \phi\subset \{|\xi|\in [2^{-2}, 2^{2}]\}$.

By repeating essentially the same argument as in Taylor \cite[(5.6), page 112]{Tay},
 we can obtain 
\beeq\label{eq-wl7}
|S_j D^s F(u)(x)|\les 2^{js} \sum_{k\in\Z}\min(1,2^{k-j}) (\M(S_k u)(x)\M H(x)+\M(H S_k u)(x))\, ,\eneq
where $H(x)\equiv G(u(x))$. 
By \eqref{eq-LPw} and \eqref{eq-wl7}, we know that for $(w_1)^{q}, (w_1 w_2)^{q}\in A_q$ with $q\in (1,\infty)$,
\begin{eqnarray}
&&\|w_1w_{2} D^s F(u)\|_{L^q}
\nonumber \ \les \
\|w_1w_{2} S_j D^s F(u)\|_{L^q\ell^2_j}\\
&\les &
\|w_1w_{2}  2^{js}  \min(1,2^{k-j}) (\M(S_k u)\M H+\M(H S_k u)) \|_{L^q\ell^2_j \ell_k^1} 
\nonumber\\
&\les &
\|w_1w_{2} 2^{ks}  \min(2^{(j-k)s},2^{(k-j)(1-s)}) (\M(S_k u)\M H+\M(H S_k u)) \|_{L^q\ell^2_j \ell_k^1} 
\nonumber\\
&\les &\|w_1w_{2} 2^{ks}  (\M(S_k u)\M H+\M(H S_k u)) \|_{L^q\ell^2_k } 
\, .\nonumber
\end{eqnarray}
For the last inequality we use Young's inequality with the assumption $s\in (0,1)$.

By  applying Minkowski's and H\"older's inequalities  to the last expression we have 
\begin{align*}
\| w_1w_2 & D^s (F( u))\|_{L^q}   \lesssim \| w_1w_2  2^{ks} \M(S_k u) \M H\|_{L^q \ell^2_k}+\|w_1w_22^{ks} \M(HS_k u) \big)\|_{L^q \ell^2_k}\\
   &\lesssim 
   \|w_2 \M H\|_{L^{q_2}}
\| w_1 2^{ks}  \M(S_k u)\|_{L^{q_1}\ell^2_k }
+
\|w_1w_{2} 2^{ks}
H S_k u \|_{L^q\ell^2_k },
\end{align*}
for any
$q_{1}, q_{2}\in (1,\infty]$ with $\frac 1q=\frac{1}{q_{1}}+\frac{1}{q_{2}}$.
For the last term in the above we used the weighted  vector valued inequality for the Hardy-Littlewood inequality which is due to Andersen-John
\cite[Theorem 3.1]{AJ80}: 
\beeq\label{M4}\int \|\M f_{j}\|_{\ell_{j}^{p}}^q wdx\le C\int \|f_{j}\|_{\ell_{j}^{p}}^q w dx, \  p,  \  q\in (1,\infty), \  w\in A_q\, .\eneq
Using \eqref{M4} again and  H\"older's inequality,
$$\| w_1w_2 D^s (F( u))\|_{L^q}   \lesssim 
   \|w_2 \M H\|_{L^{q_2}}
\| w_1 2^{ks}  S_k u\|_{L^{q_1}\ell^2_k }
+   \|w_2 H\|_{L^{q_2}}
\| w_1 2^{ks}  S_k u\|_{L^{q_1}\ell^2_k }\ .
$$
Then, by 
\eqref{eq-LPw} and its variant
$$
\| w_1 2^{ks} S_k u\|_{L^{q_1}\ell^2_k }  \equiv
\| w_1  \tilde S_k D^s u\|_{L^{q_1}\ell^2_k }\les
\| w_1  D^s u\|_{L^{q_1}}, \  \ w_1\in A_{q_1},  \ q_1\in (1,\infty),
$$
we get that
\beeq\label{eq-512-1}
\| w_1w_2 D^s (F( u))\|_{L^q}   \lesssim 
(  \|w_2 \M H\|_{L^{q_2}}
+
   \|w_2 H\|_{L^{q_2}})
\| w_1 D^{s} u\|_{L^{q_1}}.
\eneq

If  $q_2<\infty$, \eqref{eq-wLeib1} follows directly from
\eqref{eq-512-1}, by applying \eqref{M4} for the term involving $\M H$.
To handle the remaining case $q_2=\infty$, we observe that for $w\in A_1$,
\beeq\label{eq-wl9}\|w^{-1} \M H\|_{L^\infty}\les 
\|w^{-1} H\|_{L^\infty}\, .\eneq
This is trivial, since by the definition of $A_1$, we know that for a.e. $x\in \R^n$,
\begin{align*}
H(x)&\le w(x) \|w^{-1} H\|_{L^\infty}
\Rightarrow  \\ \M H(x)\le \M w(x) &\|w^{-1} H\|_{L^\infty}
\le C w(x) \|w^{-1} H\|_{L^\infty} ,
\end{align*} which gives \eqref{eq-wl9}. Together with \eqref{eq-512-1}, we get  \eqref{eq-wLeib} and this completes the proof.

\section{
Local well-posedness
}\label{Sec:3}
In this section, 
 we prove local well-posedness for
\eqref{eq-nlw} and \eqref{a-pde}, 
Theorem \ref{thm-main0} and Theorem \ref{thm-main1},
based on
Lemmas \ref{thm-trace}-\ref{thm-LE0} and Theorem \ref{thm-wLeib}.
As \eqref{eq-nlw} is a special case of \eqref{a-pde}, with $n=3$ and $p=2$, we need only
to prove Theorem \ref{thm-main1}.

\subsection{Local existence}
As usual, we prove the existence of solutions for \eqref{a-pde} through iteration. 
For fixed 
$s\in (s_{o}, 2)$ and $(u_{0}, u_{1})\in H^{s}_{\mathrm{rad}}\times H^{s-1}_{\mathrm{rad}}$, we define the iteration map
\beeq\label{eq-Picard1}
\Phi[u]:=H[u_{0}, u_{1}]+I[N[u]]\,,
\eneq where $H[\phi,\psi]=\cos (tD) \phi+
D^{-1}\sin (t
D)\psi$ is the solution map of the linear homogeneous Cauchy problem with data $(\phi,\psi)$,
$I[F](t,\cdot)=\int_0^t
D^{-1}\sin ((t-\tau)D)F(\tau,\cdot)d\tau$ is the solution map of the linear inhomogeneous Cauchy problem $\Box u=F$ with vanishing data, 
and the nonlinear term \beeq\label{eq-NL1}
N[u]:=a|\pa_t u|^p+b|\nabla_x u|^p\, .
\eneq Notice that $\Phi$ preserves radial property and 
we have $$D^{\nu}\Phi[u]=H[D^{\nu}u_{0}, D^{\nu}u_{1}]+I[D^{\nu} N[u]]\, ,$$
for any $\nu\in (0,1)$.

Let $w(x)=r^{-1+2\de}\<r\>^{-2\de-2\de_1}$  with $\de\in (0, 1/2)$ and $\de_1\in [-\de, (n-1)/2)$
 be defined by \eqref{eq-de}.
Applying Lemma \ref{thm-LE0} for $D^{\nu}\Phi[u]$ with $\nu\in [0, 1)$, we get, for $T>0$ to be determined later,
$$
\|\Phi[u]\|_{X_\nu^T} \les  \|\pa  u(0)\|_{\dot H^\nu}
+\tilde A_T\|
w^{-\frac 1 2}
 D^\nu N[u]\|_{L^2_T L^2_{x}} \ ,
$$
where
$$
\|u\|_{X_\nu^T}=
\tilde A_T^{-1}
\|w^{\frac 1 2}
 \pa D^\nu u\|_{L^2_T L^2_{x}}
+\| \pa  u\|_{L^\infty_T \dot H^\nu_{x}}\ ,
$$
$$\tilde A_T=\left\{
\begin{array}{ ll}
B_T=\min(T^\de, 1)\ ,& p>p_c\ ,\\
A_T=\min(T^\de, \sqrt{\ln(2+T)}) \ ,     &   p= p_c\ ,\\
T^\de \ ,     &   p< p_c\ . 
\end{array}
\right.
$$
As $\de\in (0, 1/2]$ and $\de_1\in [-\de, (n-1)/2)$,  we can apply Lemma \ref{thm-Ap} and Corollary \ref{thm-wLeib2} to conclude
$$\|\Phi[u]\|_{X_\nu^T}
\les  \|\pa  u(0)\|_{\dot H^\nu}
+
\tilde A_T^2
\|u\|_{X_\nu^T} 
\|w^{-1} |\pa u|^{p-1}\|_{L^{\infty}_T L^\infty_{x}}\ , \nu\in (0,1)
\, .$$
The case $\nu=0$ is also admissible, by H\"{o}lder's inequality.

We now note that the choice of $\de, \de_1$ ensures that
$$
\frac{1-2\de}{p-1}=\left\{
\begin{array}{ ll}
\frac n2-(s-1) \ ,     &   p\ge p_c\ , \\
\frac{n-1}2 \ ,     &   p< p_c\ ,
\end{array}
\right.
$$
$\frac{1+2\de_1}{p-1}\in (\frac n2-(s-1),
\frac{n-1}2]$ when
$p> p_c$ and
$\frac{1+2\de_1}{p-1}=
\frac{n-1}2$ when $p\le p_c$.
Thus
by \eqref{eq-trace1}, we see that for radial $u$,
\beeq\label{eq-trace-ap}
\|r^{\frac{1-2\de}{p-1}} \pa u\|_{L^{\infty}_T L^\infty_{x}}\les 
\left\{
\begin{array}{ ll}
\| \pa u\|_{L^{\infty}_{T} \dot H^{s-1}} \ ,     &   p\ge p_c\ , \\
\|\pa u\|_{L^{\infty}_{T}\dot H^{s-1}}^{\frac 1 2}
\| \pa u\|_{L^{\infty}_{T}\dot H^{2-s}}^{\frac 1 2}
 \ ,     &   p< p_c\ ,
\end{array}
\right.
\eneq
\beeq\label{eq-trace-ap2}
\|r^{\frac{1+2\de_1}{p-1}} \pa u\|_{L^{\infty}_T L^\infty_{x}}\les 
\left\{
\begin{array}{ ll}
\| \pa u\|_{L^{\infty}_{T}  H^{s-1}} \ ,     &   p> p_c\ , \\
\|\pa u\|_{L^{\infty}_{T}\dot H^{s-1}}^{\frac 1 2}
\| \pa u\|_{L^{\infty}_{T}\dot H^{2-s}}^{\frac 1 2}
 \ ,     &   p\le p_c\ .
\end{array}
\right.
\eneq
In conclusion, we arrive at, for any $\nu\in [0,1)$,
\beeq\label{Gl1}
\begin{split}
\|\Phi[u]\|&{}_{X_\nu^T}
\le C_{\nu}
\|\pa u(0)\|_{\dot H^\nu}      \\
&+C_{\nu} \tilde A_T^2
\|u\|_{X_\nu^T}\times
\left\{
\begin{array}{ ll}
(\|u\|_{X_{s-1}^T}+\|u\|_{X_{0}^T})^{p-1}  \ ,     &   p> p_c\ , \\
(\|u\|_{X_{s-1}^T}^2+\|u\|_{X_{s-1}^T}\|u\|_{X_{2-s}^T})^{\frac{p-1}2}
 \ ,     &   p= p_c\ ,\\
( \|u\|_{X_{s-1}^T}\|u\|_{X_{2-s}^T})^{\frac{p-1}2}
 \ ,     &   p< p_c\ .
\end{array}
\right.   
\end{split}\eneq
Similarly, 
as  $|N[u]-N[v]|\les(|\pa u|^{p-1}+|\pa v|^{p-1}) |\pa (u- v)|$,
we have
\beeq\label{Gl2}
\begin{split}
\|\Phi[u]&-\Phi[v]\|_{X_0^T}
\les\tilde A_T^2
\|u-v\|_{X_0^T}
\\
& \times
\left\{
\begin{array}{ ll}
(\|(u,v)\|_{X_{s-1}^T}+\|(u,v)\|_{X_{0}^T})^{p-1}  \ ,     &   p> p_c\ , \\
(\|(u,v)\|_{X_{s-1}^T}^2+\|(u,v)\|_{X_{s-1}^T}\|(u,v)\|_{X_{2-s}^T})^{\frac{p-1}2}
 \ ,     &   p= p_c\ ,\\
( \|(u,v)\|_{X_{s-1}^T}\|(u,v)\|_{X_{2-s}^T})^{\frac{p-1}2}
 \ ,     &   p< p_c\ .
\end{array}
\right.
\end{split}\eneq

Based on \eqref{Gl1}
with $\nu= s-1, 2-s, 0$ and
\eqref{Gl2}, it is standard to conclude that 
there exist $c, \ep_0$ such that
for $T_\ep$ given in \eqref{life}, 
$\Phi$ is a contraction mapping in the complete metric space
\[
\begin{split}
\{u\in 
L^\infty([0, T); \Hr^s)&\cap Lip( [0, T); \Hr^{s-1}),
\\
& \|\pa u\|_{L^{\infty}_{T}\dot H^{\nu}} \le 2C_\nu \|\pa u(0)\|_{\dot H^\nu}, \nu=s-1, 2-s, 0
\}
\end{split}
\]
with the metric $d(u,v)=\|u-v\|_{X_0^T}$,
and the unique fixed point is the solution we are seeking.
This completes the proof of the existence part of Theorem \ref{thm-main1}.

\subsection{Unconditional uniqueness}
Suppose we have two radial solutions $u, v\in L^\infty([0, T); \Hr^s)\cap Lip( [0, T); \Hr^{s-1})$ to \eqref{a-pde}.

Observing that $s>s_o\ge 3/2$, we have
$s-1>1/2>1/2-\de$, and so
$$\|w^{\frac 12} \pa u(t)\|_{L^2}
\les\|r^{-1/2+\de}\pa u(t)\|_{L^2}
%+\|\pa u(t)\|_{L^2}
\les \|\pa u(t)\|_{\dot H^{1/2-\de}}
\les
\|\pa u(t)\|_{H^{s-1}}
,
$$
where we have used Hardy's inequality in the second inequality. Thus,
$w^{\frac 12}\pa (u,v)\in L^2_{loc}([0,T); L^2)$.
%Let $U=u-v$, 
Then, as in the proof of the existence, we get from
 Lemma \ref{thm-LE0}, \eqref{eq-trace-ap} and \eqref{eq-trace-ap2} that, for any $T_1\in (0, T)$ and $T_{1}<1$,
\begin{align*} 
\|w^{\frac 1 2}\pa (u-v)\|_{L^2_{T_1} L^2_x}\les
 {T_1}^{2\de}\|w^{-\frac 1 2} \Box (u-v)\|_{L^2_{T_1} L^2_x} \\
 \les  {T_1}^{2\de}
 \|\pa (u,v)\|_{L^\infty_{T} H^{s-1}}^{p-1}
\|w^{\frac 12} \pa(u-v)\|_{L^2_{T_1} L^2_x}\ .
 \end{align*}
 Letting $T_1$ sufficiently small such that
 $ {T_1}^{2\de}
 \|\pa (u,v)\|_{L^\infty_{T} H^{s-1}}^{p-1}\ll 1$, we conclude that $u\equiv v$ in $[0, T_1]\times \R^n$, which yields unconditional uniqueness.
 
\subsection{Regularity} In this subsection, we prove the improved regularity of the solution $u\in C([0, T); H^s)\cap C^1( [0, T); H^{s-1})$, for which we need only to prove the continuity at $t=0$ of
$\pt^j u(t)\in \dot H^{s-j}(\R^n)$, $j=0,1$.
Applying Corollary \ref{thm-wLeib2}, \eqref{eq-trace-ap} and \eqref{eq-trace-ap2}, and
recalling the fact that $w^{\frac 1 2} \pa D^{s-1}u\in L^2_{T} L^2_x$, we have for any $T_1\in (0,T)$ and $T_{1}<1$,
    \begin{eqnarray*}
 \|\pa  u(T_1)-\pa u(0)\|_{\dot H^{s-1}_x}
          &=& \|\pa \Phi[u](T_1)-\pa u(0)\|_{\dot H^{s-1}_x}\\
&\le&\|\pa I[N[u]](T_1)\|_{\dot H^{s-1}_x}+\|\pa H[u_0, u_1](T_1)-\pa u(0))\|_{\dot H^{s-1}_x}
\\
&\les &
T_1^{\de}\|w^{-\frac 12}D^{s-1}N[u]\|_{L^2_{T_1} L^2_x}
+o(1)
\\
&\les &
  {T_1}^{\de}
 \|\pa u\|_{L^\infty_{T_1} H^{s-1}}^{p-1}
\|w^{\frac 1 2}\pa  D^{s-1}u\|_{L^2_{T_1} L^2_x}+o(1)=o(1)
  \end{eqnarray*} as $T_1\rightarrow 0+$,
which proves the continuity at $t=0$.

\section{Small data global existence with critical regularity}
In this section,  focusing on  spatial dimension two, we use the generalized Strichartz estimates to prove
Theorem \ref{thm-1}, which yields  the 
small data global existence with critical regularity.

Firstly,  we note that  $s_c=2-1/(p-1)\in (1,2)$ when  $p>5$.  Applying \eqref{eq-Stri} of Lemma \ref{thm-SSW}  and energy estimates we  get
\begin{eqnarray*}
 &  & 
 \|(\Phi[u],\pt \Phi[u])\|_{L_t^\infty (\dot H^{s_c}\times  \dot H^{s_c-1})}+\|\pa \Phi[u]\|_{L_t^{p-1} L^\infty_x} \nonumber\\
 & \les & 
\|(u_0,u_1)\|_{\dot H^{s_c}\times \dot H^{s_c-1}}+ \|N[u]\|_{L^1_t \dot H^{s_c-1}}
 \nonumber
\\
 & \les & 
 \ep 
 + \|\pa u\|_{L^{p-1}_t L^\infty_x}^{p-1} \|\pa u\|_{L^\infty_t \dot H^{s_c-1}}. 
\end{eqnarray*}
In the last inequality, we use  
Theorem \ref{thm-wLeib} 
 with $w_1=w_2=1$, $F(v)=|v|^p$, $G(v)=F'(v)$ and $v=\pa u$.
Similarly, 
by energy estimates,
we have
\begin{eqnarray*}
 \|(\Phi[u],\pt \Phi[u])\|_{L_t^\infty (\dot H^{1}\times  L^2)}
 & \les & 
\|(u_0,u_1)\|_{\dot H^{1}\times L^2}+\|N[u]\|_{L^1_t L^2_x}
 \nonumber
\\
 & \les & 
\|(u_0,u_1)\|_{\dot H^{1}\times L^2}
 + \|\pa u\|_{L^{p-1}_t L^\infty_x}^{p-1} \|\pa u\|_{L^\infty_t L^2_x} .
 \label{eq-2}
\end{eqnarray*}
Moreover, 
as $\Phi[u]$ and $\Phi[v]$ have the same data, we have
\begin{eqnarray*}
 &  & 
 \|(\Phi[u]-\Phi[v],\pt (\Phi[u]-\Phi[v]))\|_{L_t^\infty (\dot H^{1}\times  L^2)}\nonumber\\
 & \les & 
 \|N[u]-N[v]\|_{L^1_t L^2_x}
 \nonumber
\\
 & \les & 
( \|\pa u\|_{L^{p-1}_t L^\infty_x}^{p-1} +\|\pa v\|_{L^{p-1}_t L^\infty_x}^{p-1} )\|\pa (u-v)\|_{L^\infty_t L^2_x} .
 \label{eq-3}
\end{eqnarray*}
Once we have the above three estimates, it is easy and rather standard to prove the existence and uniqueness for \eqref{a-pde} in $C_t H^{s_c}\cap C_t^1 H^{s_c-1}$ with $\pa u\in L^{p-1}_t L^\infty_x$, when $(\phi,\psi)\in H^{s_c}\times H^{s_c-1}$ and $\ep$ is small enough.

When $p>3$ and the initial data are radial,  we use \eqref{eq-SSW2} of  Lemma \ref{thm-SSW} instead of \eqref{eq-Stri}. Then the same proof as before  applies without modification. This completes the proof of Theorem \ref{thm-1}.

\subsection*{Acknowledgment}
K. Hidano was supported in part by the Grant-in-Aid for
Scientific Research (C) (No. 15K04955),
Japan Society for the Promotion of Science (JSPS).
J.-C. Jiang was supported by National Sci-Tech Grant MOST 104-2115-M-007-002.
S. Lee was supported by the grant No. NRF-2015R1A4A1041675. C. Wang was supported in part by NSFC 11301478 and National Support Program for Young Top-Notch Talents.

%\bibliography{CWang815.bib}

\begin{thebibliography}{10}

\bibitem{AJ80}
K.~F. Andersen and R.~T. John.
\newblock Weighted inequalities for vector-valued maximal functions and
  singular integrals.
\newblock {\em Studia Math.}, 69(1):19--31, 1980/81.

\bibitem{CW91}
F.~M. Christ and M.~I. Weinstein.
\newblock Dispersion of small amplitude solutions of the generalized
  {K}orteweg-de {V}ries equation.
\newblock {\em J. Funct. Anal.}, 100(1):87--109, 1991.

\bibitem{FW1}
D. Fang and C. Wang.
\newblock Local well-posedness and ill-posedness on the equation of type
  {$\square u=u^k(\partial u)^{\alpha}$}.
\newblock {\em Chinese Ann. Math. Ser. B}, 26(3):361--378, 2005.

\bibitem{FW2}
D. Fang and C. Wang.
\newblock Some remarks on {S}trichartz estimates for homogeneous wave equation.
\newblock {\em Nonlinear Anal.}, 65(3):697--706, 2006.

%\bibitem{FW08}Daoyuan Fang and Chengbo Wang.\newblock Ill-posedness for semilinear wave equations with very low regularity.\newblock {\em Math. Z.}, 259(2):343--353, 2008.

\bibitem{FW11}
D. Fang and C. Wang.
\newblock Weighted {S}trichartz estimates with angular regularity and their
  applications.
\newblock {\em Forum Math.}, 23(1):181--205, 2011.

\bibitem{FW13}
D. Fang and C. Wang.
\newblock Almost global existence for some semilinear wave equations with
  almost critical regularity.
\newblock {\em Comm. Partial Differential Equations}, 38(9):1467--1491, 2013.

\bibitem{Gra14}
L. Grafakos.
\newblock {\em Classical {F}ourier analysis}, volume 249 of {\em Graduate Texts
  in Mathematics}.
\newblock Springer, New York, third edition, 2014.

\bibitem{HidTsu95}
K. Hidano and K. Tsutaya.
\newblock Global existence and asymptotic behavior of solutions for nonlinear
  wave equations.
\newblock {\em Indiana Univ. Math. J.}, 44(4):1273--1305, 1995.

\bibitem{HW17}
K. Hidano and C. Wang
\newblock 
Fractional derivatives of composite functions and the Cauchy problem for the nonlinear half wave equation.
\arXiv{1707.08319}.
% [https://arxiv.org/abs/1707.08319],



\bibitem{HWY2}
K. Hidano, C. Wang, and K. Yokoyama.
\newblock The {G}lassey conjecture with radially symmetric data.
\newblock {\em J. Math. Pures Appl. (9)}, 98(5):518--541, 2012.

\bibitem{HWY1}
K. Hidano, C. Wang, and K. Yokoyama.
\newblock On almost global existence and local well posedness for some 3-{D}
  quasi-linear wave equations.
\newblock {\em Adv. Differential Equations}, 17(3-4):267--306, 2012.

\bibitem{HiYo06}
K. Hidano and K. Yokoyama.
\newblock Space-time {$L^2$}-estimates and life span of the
  {K}lainerman-{M}achedon radial solutions to some semi-linear wave equations.
\newblock {\em Differential Integral Equations}, 19(9):961--980, 2006.


\bibitem{HZ17}
K. Hidano and D. Zha,
\newblock
Space-time $L^2$ estimates, regularity and almost global existence for elastic waves.
 \arXiv{1710.05180}.
 %[https://arxiv.org/abs/1710.05180].

\bibitem{JWY12}
J.-C. Jiang, C. Wang, and X. Yu.
\newblock Generalized and weighted {S}trichartz estimates.
\newblock {\em Commun. Pure Appl. Anal.}, 11(5):1723--1752, 2012.

\bibitem{JK84}
F.~John and S.~Klainerman.
\newblock Almost global existence to nonlinear wave equations in three space
  dimensions.
\newblock {\em Comm. Pure Appl. Math.}, 37(4):443--455, 1984.

\bibitem{Jo81}
F. John.
\newblock Blow-up for quasilinear wave equations in three space dimensions.
\newblock {\em Comm. Pure Appl. Math.}, 34(1):29--51, 1981.

\bibitem{Ka95}
T. Kato.
\newblock On nonlinear {S}chr\"odinger equations. {II}. {$H^s$}-solutions and
  unconditional well-posedness.
\newblock {\em J. Anal. Math.}, 67:281--306, 1995.

\bibitem{KSS}
M. Keel, H.~F. Smith, and C.~D. Sogge.
\newblock Almost global existence for some semilinear wave equations.
\newblock {\em J. Anal. Math.}, 87:265--279, 2002.
\newblock Dedicated to the memory of Thomas H. Wolff.


\bibitem{KlMa93}
S. Klainerman and M. Machedon.
\newblock Space-time estimates for null forms and the local existence theorem.
\newblock {\em Comm. Pure Appl. Math.}, 46(9):1221--1268, 1993.


\bibitem{KM95}
S.~Klainerman and M.~Machedon.
\newblock Smoothing estimates for null forms and applications.
\newblock {\em Duke Math. J.}, 81(1):99--133 (1996), 1995.
\newblock A celebration of John F. Nash, Jr.

\bibitem{KM96}
S. Klainerman and M. Machedon.
\newblock Estimates for null forms and the spaces {$H_{s,\delta}$}.
\newblock {\em Internat. Math. Res. Notices}, (17):853--865, 1996.

\bibitem{KS97}
S. Klainerman and S. Selberg.
\newblock Remark on the optimal regularity for equations of wave maps type.
\newblock {\em Comm. Partial Differential Equations}, 22(5-6):901--918, 1997.

\bibitem{Ld93}
H. Lindblad.
\newblock A sharp counterexample to the local existence of low-regularity
  solutions to nonlinear wave equations.
\newblock {\em Duke Math. J.}, 72(2):503--539, 1993.

\bibitem{Ld96}
H. Lindblad.
\newblock Counterexamples to local existence for semi-linear wave equations.
\newblock {\em Amer. J. Math.}, 118(1):1--16, 1996.

\bibitem{LMSTW}
H. Lindblad, J. Metcalfe, C.~D. Sogge, M. Tohaneanu, and
  C. Wang.
\newblock The {S}trauss conjecture on {K}err black hole backgrounds.
\newblock {\em Math. Ann.}, 359(3-4):637--661, 2014.

\bibitem{LW17}
M. Liu and C. Wang, 
\newblock 
Concerning ill-posedness for semilinear wave equations.
\newblock
in preparation.

\bibitem{MetSo06}
J. Metcalfe and C.~D. Sogge.
\newblock Long-time existence of quasilinear wave equations exterior to
  star-shaped obstacles via energy methods.
\newblock {\em SIAM J. Math. Anal.}, 38(1):188--209, 2006.

\bibitem{MeTa12MA}
J. Metcalfe and D. Tataru.
\newblock Global parametrices and dispersive estimates for variable coefficient
  wave equations.
\newblock {\em Math. Ann.}, 353(4):1183--1237, 2012.

\bibitem{MuSch13}
C. Muscalu and W. Schlag.
\newblock {\em Classical and multilinear harmonic analysis. {V}ol. {I}}, volume
  137 of {\em Cambridge Studies in Advanced Mathematics}.
\newblock Cambridge University Press, Cambridge, 2013.


\bibitem{Plan03}
F. Planchon.
\newblock On uniqueness for semilinear wave equations.
\newblock {\em Math. Z.}, 244(3):587--599, 2003.

\bibitem{PS}
G. Ponce and T.~C. Sideris.
\newblock Local regularity of nonlinear wave equations in three space
  dimensions.
\newblock {\em Comm. Partial Differential Equations}, 18(1-2):169--177, 1993.

\bibitem{Ram97}
M.~A. Rammaha. \newblock 
A note on a nonlinear wave equation in two and three space dimensions.\newblock {\em Comm. Partial Differential Equations},
22 (1997), no. 5-6, 799--810.


\bibitem{Si83}
T.~C. Sideris.
\newblock Global behavior of solutions to nonlinear wave equations in three
  dimensions.
\newblock {\em Comm. Partial Differential Equations}, 8(12):1291--1323, 1983.


\bibitem{SSW12}
H.~F. Smith, C.~D. Sogge, and C. Wang.
\newblock Strichartz estimates for {D}irichlet-wave equations in two dimensions
  with applications.
\newblock {\em Trans. Amer. Math. Soc.}, 364(6):3329--3347, 2012.

\bibitem{SmTa05}
H.~F. Smith and D. Tataru.
\newblock Sharp local well-posedness results for the nonlinear wave equation.
\newblock {\em Ann. of Math. (2)}, 162(1):291--366, 2005.

\bibitem{Sta97}
G. Staffilani.
\newblock On the generalized {K}orteweg-de {V}ries-type equations.
\newblock {\em Differential Integral Equations}, 10(4):777--796, 1997.

\bibitem{Ster05}
J. Sterbenz.
\newblock Angular regularity and {S}trichartz estimates for the wave equation.
\newblock {\em Int. Math. Res. Not.}, (4):187--231, 2005.
\newblock With an appendix by Igor Rodnianski.

\bibitem{Ster07}
J. Sterbenz.
\newblock Global regularity and scattering for general non-linear wave
  equations. {II}. {$(4+1)$} dimensional {Y}ang-{M}ills equations in the
  {L}orentz gauge.
\newblock {\em Amer. J. Math.}, 129(3):611--664, 2007.

\bibitem{Tataru99}
D. Tataru.
\newblock On the equation {$\square u=|\nabla u|^2$} in {$5+1$} dimensions.
\newblock {\em Math. Res. Lett.}, 6(5-6):469--485, 1999.

\bibitem{Tay}
M.~E. Taylor.
\newblock {\em Tools for {PDE}}, volume~81 of {\em Mathematical Surveys and
  Monographs}.
\newblock American Mathematical Society, Providence, RI, 2000.
\newblock Pseudodifferential operators, paradifferential operators, and layer
  potentials.

\bibitem{Tz98}
N. Tzvetkov.
\newblock Existence of global solutions to nonlinear massless {D}irac system
  and wave equation with small data.
\newblock {\em Tsukuba J. Math.}, 22(1):193--211, 1998.

\bibitem{W15}
C. Wang.
\newblock The {G}lassey conjecture on asymptotically flat manifolds.
\newblock {\em Trans. Amer. Math. Soc.}, 367(10):7429--7451, 2015.

\bibitem{Zh00}
Y. Zhou.
\newblock Uniqueness of generalized solutions to nonlinear wave equations.
\newblock {\em Amer. J. Math.}, 122(5):939--965, 2000.


\bibitem{Zh01}
Y. Zhou.
\newblock Blow up of solutions to the {C}auchy problem for nonlinear wave
  equations.
\newblock {\em Chinese Ann. Math. Ser. B}, 22(3):275--280, 2001.

\end{thebibliography}

\end{document}